\documentclass[letterpaper]{article}

\usepackage{amsmath, amssymb}
\usepackage{amsthm} 
\theoremstyle{plain} 
\newtheorem{theorem}{\indent\sc Theorem}[section] 
\newtheorem{Lemma}[theorem]{\indent\sc Lemma}
\newtheorem{Colr}[theorem]{\indent\sc Corollary}
\newtheorem{prop}[theorem]{\indent\sc Proposition}

\theoremstyle{definition}


\title{Deformation equivalence classes of complex surfaces with the first Betti number one, and the second Betti number zero} 
\author{Murakami Shota} 

\usepackage{latexsym}
\usepackage{graphics}
\usepackage{multicol}
\usepackage{amscd}
\usepackage{amsmath,amsthm,amssymb,euscript}
\usepackage{ascmac}
\usepackage{tabularx,array}
\usepackage{bm}
\usepackage{hhline}
\usepackage{amsmath,amscd}
\usepackage[dvips]{color} 
\usepackage{enumerate}

\def\Z{\mathbb{Z}}

\def\R{\mathbb{R}}
\def\C{\mathbb{C}}

\def\HC{\mathbb{H}}

\def\Aut{{\rm Aut} }
\def\GMinus{G^{-}_{N,p,q,r}}
\def\Gminus2{G^{-}_{N',p',q',r'}}
\def\Endm{{\rm End} }
\def\nbsp{\hspace{1cm}}

\newcommand{\CP}{\mathbb{C}\mathbb{P}}

\makeatother

\date{}
\begin{document}
\maketitle
\footnote{ 
2010 \textit{Mathematics Subject Classification}.
Primary 32J15; Secondary 14J25.
}
\footnote{ 
\textit{Key words and phrases}. 
Deformation equivalence, Inoue surface
}
\begin{abstract}
We will prove that the number of deformation equivalence classes of surfaces homotopy equivalent to a smooth, closed 4-manifold is finite, if the first Betti number is equal to one, and the second Betti number is equal to zero.
\end{abstract}

\section*{Introduction}
\indent \indent Recently, the study of compact complex surfaces has made many outstanding progresses. One of the topics in this area is to investigate the discrepancy between homotopy class, homeomorphism class, diffeomorphism class and deformation equivalence classes of surfaces.  In this paper we ask whether the number of deformation equivalence classes of surfaces  homotopy equivalent (or diffeomorphic) to a smooth closed 4-manifold is finite or not. For this problem, the celebrated Yau's result \cite{Yau} states that any surface homotopy equivalent to $\CP^2$ is biholomorphic to $\CP^2$, hence there exists only one deformation equivalence class of surfaces homotopy equivalent to $\CP^2$. According to \cite{Friedman},  the number of deformation equivalence classes of surface diffeomorphic to $M$ is finite, if $M$ is a  smooth, closed 4-manifold $M$ with $b_1(M) \neq 1$. It is also known that the number of deformation equivalence classes of surfaces homotopy equivalent to $M$ is finite if and only if  $M$ is not homotopy equivalent to an elliptic surface whose fundamental group is finite cyclic. 
Our result is an affirmative answer for $M$ with $b_1(M) = 1$ and $b_2(M) = 0$. 

\paragraph{Main Theorem} Let $M$ be a smooth, closed 4-manifold with $b_1 (M) =1$ and $b_2 (M) = 0$. Then the number of deformation equivalence classes of surfaces homotopy equivalent to $M$ is finite. \\

\indent  Our strategy is as follows. Let $S$ be a compact surface with $b_1 (S) =1$ and $b_2 (S) = 0$.  By the Enriques-Kodaira classification \cite{VanDeVan} and the results of Bogomolv \cite{Bogomolov}, and Teleman \cite{Teleman}, $S$ is biholomorphic to an elliptic surface, a Hopf surface, or an Inoue surface. We first assume that $S$ is an Inoue surface. We show that $S$ cannot deform to surfaces other than Inoue surfaces by comparing the fundamental groups. The fundamental group also distinguishes all the Inoue surfaces homotopy equivalent to $S$. We show that the number of deformation equivalence classes is finite by constructing biholomorphic maps explicitly. Then, we will prove the theorem for elliptic surfaces and Hopf surfaces using the results of \cite{Friedman}.\\
\indent We remark that the classification of surfaces with $b_1 =1 , b_2 >0$ and Kodaira dimension $-\infty$, remains open. It is conjectured  that $S$ contains a global spherical shell \cite{Teleman_arxiv} if $S$ is minimal. If this conjecture is proved affirmatively,  \cite{Friedman} and \cite{Kato} implies that any surface with $b_1 =1$ and $b_2 >0$ is deformation equivalent to a surface with $b_1 =1$, $b_2 = 0$ blown up finitely times. Hence, by our main theorem, it follows that the number of deformation equivalence classes of surfaces diffeomorphic to a closed 4-manifold is finite. 
\section{Inoue surfaces of type $S^0$}
\indent \indent In \cite{Inoue_Original}, Inoue constructed surfaces with $b_1 =1, b_2 =0$. We follow \cite{Hasegawa_Solvmanifold} to divide them into three types; type  $S^0$,$S^{+}$, and $S^{-}$. We begin by recalling the definition of Inoue surfaces of type $S^0$.\\
\indent Let $M = (m_{ij} ) \in SL(3, \Z) $ be a matrix with eigenvalues $\alpha>1 , \beta , \bar{\beta} \in \C \setminus \R$.  We choose eigenvectors $\bm{a} = {}^{T} (a_1 , a_2, a_3) \in \R^3$ and $\bm{b} = {}^{T}(b_1 ,b_2,b_3) \in \C^3$ of $M$ corresponding to $\alpha$ and $\beta$ respectively. Define $G_M$ to be the subgroup of $\Aut(\HC \times \C)$ generated by
\begin{eqnarray}
g_0 & : & (w,z) \mapsto ( \alpha w,\beta z) \label{eqn:SMDef1}, \\
g_i  & : & (w,z) \mapsto (w + a_i, z + b_i), \label{eqn:SMDef2}
\end{eqnarray}
where $i=1,2,3,$ and $\HC$ is the upper half complex plane. The quotient surface $S^0_{M, \bm{a},\bm{b}} = \HC \times \C / G_M$ is called an Inoue surface of type $S^0$. \\
\indent Recall that two complex manifolds $S_1$ and $S_2$ are deformation equivalent if there exists connected complex manifolds $\chi$ and $B$, a smooth proper holomorphic map $\pi:\chi \to B$, and two points $t_1,t_2 \in B$ such that $\pi^{-1}(t_i)$ is biholomorphic to $S_i$. In \cite{Inoue_SM}, Inoue showed that if $S$ is a surface whose fundamental group is isomorphic to $\pi_1(S^0_{M, \bm{a},\bm{b}})$, then $S$ is biholomorphic either to $S^0_{M, \bm{a},\bm{b}}$ or $S^{0}_{M,\bm{a},\bar{\bm{b}}}$. He also showed that $S^0_{M, \bm{a},\bm{b}}$ and  $S^{0}_{M,\bm{a},\bar{\bm{b}}}$ are not deformation equivalent. Hence, the number of deformation equivalence classes of surfaces homotopy equivalent to $S$ is exactly two.

\section{Inoue surfaces of type $S^+$}
\indent \indent Let $N = (n_{ij} ) \in SL(2, \Z) $ be a matrix with eigenvalues $\alpha>1, 1/ \alpha$. We choose eigenvectors $\bm{a} = {}^{T} (a_1, a_2) \in \R^2$ and $\bm{b} = {}^{T}(b_1,b_2) \in \R^2$ of $N$ corresponding to $\alpha$ and $1/\alpha$ respectively. Fix $t \in \C$ and integers $p,q, r,$ where $r \neq 0$. Let $\theta = \det( \bm{a}, \bm{b} )$. Define $c_1, c_2 \in \R$ to be solution of the following equation
\begin{eqnarray*}
(N - I) \left( \begin{array}{c} c_1 \\ c_2 \end{array} \right) +  \left( \begin{array}{c} e_1 \\ e_2 \end{array} \right) - \frac{\theta}{r}  \left( \begin{array}{c} p \\ q \end{array} \right) = 0,
\end{eqnarray*}
where $e_i = \frac{1}{2}n_{i1} (n_{i1} -1) a_1 b_1 + \frac{1}{2}n_{i2} (n_{i2} -1) a_2 b_2 + n_{i1} n_{i2} b_1 a_2$, and $i = 1,2$. Let 
$G^{+}_{N,p,q,r}$ be the subgroup of $\Aut(\HC \times \C)$ generated by
\begin{eqnarray*}
g_0: (w,z) & \mapsto & (\alpha w, z + t), \\
g_i : (w,z) & \mapsto & (z + a_i, w + b_i z + c_i ), \\
g_3: (w,z) & \mapsto & (z, w - \theta /r),
\end{eqnarray*}
where $i=1,2$. The quotient surface $S^{+}_{N,p,q,r,t,\bm{a},\bm{b}} = \HC \times \C / G^{+}_{N,p,q,r}$ is called an Inoue surface of type  $S^{+}$. \\
\indent We have the following relations:
\begin{eqnarray}
g_0 g_1 g_0^{-1} &=& g_1^{n_{11}} g_2^{n_{12}} g_3^p, \label{eqn:SPlus_Rel1} \\ 
g_0 g_2 g_0^{-1} &=& g_1^{n_{21}} g_2^{n_{22}} g_3^q, \label{eqn:SPlus_Rel2} \\
 g_1 g_2 g_1^{-1} g_2^{-1} &=& g_3^r  \label{eqn:SPlus_Rel3}, \\
 g_i g_3 &=& g_3 g_i \label{eqn:SPlus_Rel4},
\end{eqnarray}
where $i=0,1,2,3$. Any element $g \in G^{+}_{N,p,q,r}$ can be written uniquely as a product $g = \prod_{i=0}^{3} g_i^{n_i}$, where $n_0,n_1,n_2,n_3$ are integers. Thus, $G^{+}_{N,p,q,r}$  has a presentation in terms of the generators $g_0,g_1,g_2,g_3$, with relations (\ref{eqn:SPlus_Rel1}) - (\ref{eqn:SPlus_Rel4}). The center of $G^{+}_{N,p,q,r}$ is an infinite cyclic group generated by $g_3$. 
\paragraph{Notation} Let $G$ be a group. We denote by $\Gamma_G$ the normal subgroup
 \begin{gather}
 \label{eqn:Gamma_Group}
\{ g \in G | g  \textrm{ mod} [G,G] \textrm{ is of finite order} \},
\end{gather}
where $[G,G]$ is the commutator subgroup. \\

The subgroup $\Gamma_{G^{+}_{N,p,q,r}}$ of $G^{+}_{N,p,q,r}$ has the following properties:
\begin{itemize}
\item[(i)] $\Gamma_{G^{+}_{N,p,q,r}}$ is generated by $g_1,g_2,g_3$,
\item[(ii)] the center of $\Gamma_{G^{+}_{N,p,q,r}}$ is an infinite cyclic group generated by $g_3$,
\item[(iii)]  the quotient group of $\Gamma_{G^{+}_{N,p,q,r}}$ by its center is a free abelian group of rank two generated by the classes of $g_1$ and $g_2$,
\item[(iv)] the quotient group $G^{+}_{N,p,q,r}/\Gamma_{G^{+}_{N,p,q,r}}$ is an infinite cyclic group generated by the class of $g_0$.
\end{itemize}
\begin{Lemma}
\label{thm:SPlus_diffeo}
Inoue surfaces $S^{+}_{N,p,q,r,t,\bm{a},\bm{b}}$ and $S^{+}_{N',p',q',r',t',\bm{a'},\bm{b'}}$ constructed above are homotopy equivalent if and only if there exists $K \in GL(2,\Z) , \delta, \epsilon \in \{ 1,-1 \}$ such that
\begin{eqnarray}
r' & = & \delta \det K r, \label{eqn:SPlus_diffeo1}\\
N' & = & KN^{\epsilon} K^{-1}, \label{eqn:SPlus_diffeo2} \\
\delta \left( \begin{array}{c}
p' \\
q'
\end{array} \right) - \epsilon  K \left(
\begin{array}{c}
p \\
q
\end{array}
 \right) & \in & r \Z^2 + (N' - I) \Z^2. \label{eqn:SPlus_diffeo3}
\end{eqnarray}
\end{Lemma}

\paragraph{Proof}
 Assume that $S^{+}_{N,p,q,r,t,\bm{a},\bm{b}}$ and $S^{+}_{N',p',q',r',t',\bm{a'},\bm{b'}}$ are homotopy equivalent. Then, there exists an isomorphism $\varphi: G^{+}_{N',p',q',r'} \to G^{+}_{N,p,q,r} $. By property (ii), $\varphi ( g'_3 )  = g_3^{\delta}$ for some $\delta \in \{ 1, -1 \}$. From property (i), there exists integers $k_{i1},k_{i2},k_{i3} $  such that
\begin{eqnarray}
\varphi (g'_{i}) & = &  g_1^{k_{i1}} g_2^{k_{i2}} g_3^{k_{i3}}, \label{eqn:SPlusDiffeoGroup2}
\end{eqnarray}
where $i =1,2$. We infer from property (iii) that if we let $K = (k_{ij})_{i,j=1,2}$, then $K \in GL(2,\Z)$. Property (iv) implies that there exists $\epsilon \in \{ 1,-1 \}, $ and integers $l_1,l_2,l_3$ such that
\begin{eqnarray}
\varphi (g'_{0}) &= & g_0^{\epsilon} g_1^{l_1} g_2^{l_2} g_3^{l_3}. \label{eqn:SPlusDiffeoGroup1}
 \end{eqnarray} 

By relation (\ref{eqn:SPlus_Rel3}), we obtain (\ref{eqn:SPlus_diffeo1}). For the sake of simplicity, let $ n'_{13} = p' , n'_{23} = q'$. Suppose $\epsilon = 1$. By (\ref{eqn:SPlusDiffeoGroup2}) and (\ref{eqn:SPlusDiffeoGroup1}), we have
\begin{eqnarray*}
\varphi ( g'_0 g'_i {g'_0}^{-1} ) & = & g_1^{(k_{i1} , k_{i2} ){}^{T}(n_{11},n_{21})} g_2^{(k_{i1} , k_{i2} ){}^{T}(n_{12},n_{22})} g_3^{k_{i3} + (k_{i1} , k_{12}){}^{T}(p,q) + r  (k_{i1} , k_{i2} )  {}^{T}(-l_2, l_1) +r  a }, \\ 
\varphi ( {g'}_1^{n'_{i1}} {g'}_2^{n'_{i2}} {g'}_3^{n'_{i3}}) & = & g_1^{(n'_{i1} , n'_{i2}) {}^{T}(k_{11}, k_{21})} g_2^{(n'_{i1} , n'_{i2}) {}^{T}(k_{12}, k_{22})} g_3^{(n'_{i1} , n'_{i2}) {}^{T}(k_{13}, k_{23}) + \delta n'_{i3} + r b},
\end{eqnarray*}
where $i=1,2$ and $a,b$ are integers which depends only on $K$,$N$, and $N'$. Thus,

\begin{eqnarray*}
KNK^{-1} &=& N', \\
\delta \left( \begin{array}{c}
p' \\ q'
\end{array} \right) - K \left( \begin{array}{c}
p \\ q
\end{array} \right) & \in &  (N' - I) \Z^2 + r \Z^2.
\end{eqnarray*}
Similarly, if $\epsilon = -1$,  we obtain $KN^{-1} K^{-1} = N'$ and
\begin{eqnarray*}
\delta \left( \begin{array}{c}
p' \\ q'
\end{array} \right) + K \left( \begin{array}{c}
p \\ q
\end{array} \right) \in  (N' - I) \Z^2 + r \Z^2.
\end{eqnarray*}
Hence, we get (\ref{eqn:SPlus_diffeo2}) and (\ref{eqn:SPlus_diffeo3}). \\
\indent Next we show the converse. Recall that the diffeomorphism type of solvmanifolds is determined by the fundamental group. Since $S^{+}_{N,p,q,r,t,\bm{a},\bm{b}}$ and  $S^{+}_{N',p',q',r',t',\bm{a'},\bm{b'}}$ are both solvmanifolds \cite{Hasegawa_Solvmanifold}, it suffices to show that $G^{+}_{N',p',q',r'}$ and $G^{+}_{N,p,q,r}$ are isomorphic. By (\ref{eqn:SPlus_diffeo3}), there exists $k_{13},k_{23} , u , v \in \Z$ such that
\begin{eqnarray*}
\delta \left( \begin{array}{c}
p' \\ q'
\end{array} \right) - \epsilon K \left( \begin{array}{c}
p \\ q
\end{array} \right) =  (N' - I) \left( \begin{array}{c}
k_{13} \\ k_{23}
\end{array} \right) + r \left( \begin{array}{c}
u \\ v
\end{array} \right).
\end{eqnarray*}
Suppose $\epsilon = 1$. Let $l_1,l_2$ be integers which satisfy the following equation:
\begin{eqnarray*}
\left( \begin{array}{c}
-l_2 \\ l_1
\end{array} \right)= K^{-1} \left( \begin{array}{c}
u \\ v
\end{array} \right) - \left( \begin{array}{c}
a \\ b
\end{array} \right).
\end{eqnarray*}
The desired isomorphism $\varphi :G^{+}_{N',p',q',r'} \to G^{+}_{N,p,q,r}$ is given by
\begin{eqnarray*}
\varphi (g'_{0}) &= & g_0^{\epsilon} g_1^{l_1} g_2^{l_2}, \\
\varphi (g'_{i}) & = &  g_1^{k_{i1}} g_2^{k_{i2}} g_3^{k_{i3}}, \\
\varphi ( g'_3 )  &=& g_3^{\delta},
\end{eqnarray*}
where $i=1,2$. The case $\epsilon = -1$ can be proved similarly.
\qed \\
 
In \cite{Surjective_Isom},  Fujimoto and Nakayama introduced the group $\Gamma_r$ to  construct a surjective holomorphic map from $ S^{+}_{N,p,q,r,t,\bm{a},\bm{b}}$ to itself. We apply their method to construct a biholomorpic map between two Inoue surfaces of type $S^{+}$ with different parameters. However, we will assume that $r = r'$, which is enough to show our main theorem. We first review its definition. Let $\Gamma_r = \Z^2 \times \Z [r/2]$, where $\Z [r/2] = \Z$ if $r$ is even, and $\frac{1}{2} \Z$ if $r$ is odd.  The group law of $\Gamma_r$ is defined as follows:
\begin{eqnarray}
\label{eqn:Gamma_to_Gammar}
(\zeta, y) (\zeta',y') = ( \zeta + \zeta', y + y' +\frac{r}{2} \det( \zeta, \zeta')),
\end{eqnarray}
where $\zeta, \zeta' \in \Z^2$ are row vectors, and $y,y' \in \Z [r/2]$. By (8.4) of \cite{Surjective_Isom}, we have
\begin{eqnarray}
\label{SPlus_SPEQ}
(N - I) \bm{c} = \frac{\theta}{r} \bm{p},
\end{eqnarray}
where $\bm{c}= {}^{T} (c_1 - a_1 b_1/2,  c_2 - a_2 b_2/2),$ and $\bm{p} = {}^{T} (p + (r/2) n_{11} n_{12}, q+(r/2) n_{21} n_{22})$. Let $\mu: \Gamma_{G^{+}_{N,p,q,r}} \to \Gamma_r$ be an injective homomorphism given by
 \begin{eqnarray}
 \label{eqn:Gammar_HCAction}
g_1^{l_1} g_2^{l_2} g_3^{l_3} \mapsto ((l_1,l_2),l_3 + l_1 l_2 \frac{r}{2}).
\end{eqnarray} The group $\Gamma_r$ acts on $\HC \times \C$ as follows:
\begin{eqnarray}
\label{eqn:GammaR_Action}
(w,z) \mapsto  (\zeta,y) (w,z) = \left( w + \zeta \bm a, z + (\zeta \bm b) w + \zeta \bm {c} - \frac{\theta}{r} y + \frac{1}{2} (\zeta \bm a) (\zeta \bm b) \right).
\end{eqnarray}
The action above factors through $\Gamma_r$ by $\mu$.

%
%
%
%
%
%
%
%
%
%
%
%
%
%
%
\indent Let us review the endomorphisms of $\Gamma_r$. According to Lemma 8.4 of \cite{Surjective_Isom}, an endomorphism $\varphi$ of $\Gamma_r$ can be written as
  \begin{eqnarray*}
\varphi (\zeta, y) = (\zeta K, \zeta \bm{v} + (\det K) y),
\end{eqnarray*}
for some $K = (k_{ij}) \in M_2(\Z)$ and $\bm{v} = {}^T (v_1,v_2) \in \Z [r/2]^2$. Furthermore, $\Endm(\Gamma_r)$ is anti-isomorphic to the semi-group $M_2 (\Z) \times \Z [r/2]^2$ via the map
\begin{eqnarray}
\label{GammaRendIsom}
\Endm(\Gamma_r) \to M_2 (\Z) \times \Z [r/2]^2 : \varphi \mapsto (K,v),
\end{eqnarray}
where the multiplicative structure on $M_2 (\Z) \times \Z [r/2]^2$ is defined by 
\begin{eqnarray*}
(K_1, \bm{v_1} ) (K_2, \bm{v_2}) = (K_1 K_2, K_1 \bm{v_2} + (\det K_2) \bm{v_1}).
\end{eqnarray*}
%
%
Throughout this paper, we identify $\Endm(\Gamma_r)$ with $M_2 (\Z) \times \Z [r/2]^2$ by (\ref{GammaRendIsom}).
\begin{Lemma} %
\label{lem:Gamma_R}
 Assume that $r= r'$.
\begin{enumerate} 
\item[{\normalfont(1)}] Any homomorphism $\Gamma_{G^{+}_{N',p',q',r'}} \to \Gamma_{G^{+}_{N,p,q,r}}$ lifts to a unique endomorphism of $\Gamma_r$. 
\item[{\normalfont(2)}] The pair $(K, \bm{v}) \in M_2 (\Z) \times \Z [r/2]^2$ is induced from a homomorphism $\Gamma_{G^{+}_{N',p',q',r'}} \to \Gamma_{G^{+}_{N,p,q,r}}$ if and only if $v_1 - (r/2) k_{11} k_{12}, v_2 - (r/2) k_{21} k_{22} \in \Z$ where $K = (k_{ij} )$, and $\bm{v} = {}^T (v_1, v_2)$.
\item[{\normalfont(3)}] The lift of the automorphism $\nu: \gamma \mapsto g_0 \gamma g_0^{-1}$ of $\Gamma_{G^{+}_{N,p,q,r}}$ corresponds to $(N, \bm{p})$.
\item[{\normalfont(4)}] Let $\rho:{G^{+}_{N',p',q',r'}} \to {G^{+}_{N,p,q,r}}$ be an isomorphism such that $\rho(g'_0) = g_0 g_1^{l_1} g_2^{l_2}$ for some integers $l_1,l_2$. Denote $(K, \bm{v}) \in GL(2,\Z) \times \Z[r/2]^2$ the lift of $\rho|{\Gamma_{{G^{+}_{N',p',q',r'}}}}$. Then $(K, \bm{v})$ satisfies the following conditions:
\begin{eqnarray}
KN = N'K,  \nbsp  (N' - I ) \bm{v}  = K \bm{p} - (\det K) \bm{p}' + Kr \left( \begin{array}{c} -l_2 \\ l_1 \end{array} \right) \label{Splus_GammaR_equiv2} \label{Gamma_r_isom}, \\
v_i - (r/2) k_{i1} k_{i2}  \in \Z,  \nbsp \nbsp \nbsp \nbsp \label{Gamma_r_isom2}
\end{eqnarray}
where $\bm{v} = {}^T (v_1, v_2) $. Conversely, assume that $(K, \bm{v}) \in GL(2,\Z) \times \Z[r/2]^2$ satisfies $(\ref{Splus_GammaR_equiv2})$ and $(\ref{Gamma_r_isom2})$ for some integers $l_1$ and $l_2$. Let $\varphi:\Gamma_{G^{+}_{N',p',q',r'}} \to \Gamma_{G^{+}_{N,p,q,r}}$ be a isomorphism whose lift corresponds to $(K, \bm{v})$. Then, the map
\begin{eqnarray*}
 \rho:G^{+}_{N',p',q',r'} \to G^{+}_{N,p,q,r}: {g'_0}^{l_0} \gamma \mapsto (g_0 g_1^{l_1} g_2^{l_2})^{l_0} \varphi (\gamma),
 \end{eqnarray*}
 is an isomorphism, where $\gamma \in \Gamma_{G^{+}_{N',p',q',r'}}$.
\end{enumerate}
\end{Lemma}
%
%
%
%
%
%
%
%
%
%
%
%
%
%
%
%
%
%

\paragraph{Proof}
Let $\varphi: \Gamma_{G^{+}_{N',p',q',r'}} \to \Gamma_{G^{+}_{N,p,q,r}}$ be a homomorphism. Then $\varphi(g'_i) = g_1^{k_{i1}} g_2^{k_{i2}} g_3^{k_{i3}}$ for some integers $k_{ij}$, where $i=1,2$ and $j=1,2,3$. If we let $K = (k_{ij})_{i,j=1,2}$, we see that the endomorphism of $\Gamma_r$ corresponding to 
\begin{eqnarray*}
\left( K, \left( \begin{array}{c}  k_{13} + \frac{r}{2} k_{11} k_{12}  \\  k_{23} + \frac{r}{2}k_{21} k_{22} \end{array} \right) \right) \in GL(2,\Z) \times \Z [r/2]^2
\end{eqnarray*}
is the lift of $\varphi$. This shows (1). It is obvious that $((l_1, l_2), \lambda) \in \Gamma_r$ is in the image of $\mu $ if and only if $\lambda - (r/2) l_1 l_2 \in \Z$, hence we have (2). (3) follows from relations (\ref{eqn:SPlus_Rel1}) and (\ref{eqn:SPlus_Rel2}). For (4), let $\iota$ be the automorphism of $\Gamma_{G^{+}_{N,p,q,r}}$ given by the conjugation of $\eta = g_1^{l_1} g_2^{l_2}$. We define the automorphism $\nu'$ of $\Gamma_{G^{+}_{N',p',q',r'}}$ similarly to $\nu$.  Let $\varphi = \rho|\Gamma_{G^{+}_{N',p',q',r'}}$. Then, $\rho( g'_0 \gamma {g'}_0^{-1}) = g_0 \eta \rho ( \gamma) \eta^{-1} g_0^{-1}$ implies
\begin{eqnarray}
\varphi \circ \nu' = \nu \circ \iota \circ \varphi. \label{eqn:nuandphi}
\end{eqnarray}
Since the lift of $\iota$ corresponds to $(I, r {}^T (-l_2, l_1))$,  we obtain 
\begin{eqnarray*}
(N', \bm{p'}) (K, \bm{v}) = (K, \bm{v})  (I, r {}^T (-l_2, l_1)) (N, \bm{p}),
\end{eqnarray*} which is equivalent to (\ref{Gamma_r_isom}). \\
\indent We now show the converse.  By assumption, we have (\ref{eqn:nuandphi}). It follows that $\rho$ is an isomorphism by calculation. \qed
%
%
%
%
%
%
%
%
%
%
%
%
%
%
%
%

%
%
%
%
%
%
%
%
%
%
%
%
%
%
%
%
%
%
%
%
%
%
%
%
%
%
%
%
%
%
%
%
%
%
%
%
%
%
%

\begin{Lemma}
\label{lem:SPlusBihol}
Assume that $r=r'$ and that there exists $K \in GL(2,\Z)$,$c>0 , u,v \in \Z$ such that
\begin{gather}
KNK^{-1}  =   N' ,\label{eqn:SPlus_Bihol_2} \\
c \bm{a'}  =  K \bm{a} ,  f \bm{b'}  =  cK \bm{b}, \label{eqn:SPlus_Bihol_4}\\
ft' - t   =    \frac{- \alpha}{\alpha -1} (\eta \bm{a})(\eta \bm{b}) + \eta \bm{c} + \frac{1}{2} (\eta \bm{a} ) (\eta \bm{b}) - \frac{\theta}{r}  l_1 l_2, \label{eqn:SPlus_Bihol_6} \\
\eta \in \Z^2, \nonumber
\end{gather}
where $f = \displaystyle \frac{det K \theta}{\theta'} , $ and $\eta = (l_1 , l_2)$ is a solution of the following equation:
\begin{eqnarray}
\label{eqn:GPlus_Bihol_GammaR_diffeo}
(N' - I) \left( \begin{array}{c} u +\frac{r}{2} k_{11} k_{12} \\ v + \frac{r}{2}  k_{21} k_{22} \end{array}\right) + (\det K) \bm{p'} = K \bm{p} + Kr \left( \begin{array}{c} -l_2 \\ l_1 \end{array} \right).
\end{eqnarray}
Then  $S^{+}_{N,p,q,r,t,\bm{a},\bm{b}}$ and $S^{+}_{N',p',q',r',t',\bm{a'},\bm{b'}}$ are biholomorphic.
\end{Lemma}
\paragraph{Proof}
Let $\varphi:\HC \times \C \to \HC \times \C:(w,z) \mapsto (cw + d, ew +fz),$ where 
\begin{eqnarray}
\displaystyle d = - \frac{\alpha }{\alpha -1} (\eta \bm{a}), \displaystyle e =  \frac{c}{\alpha -1} (\eta \bm{b}). \label{eq:Splus_biHolE}
\end{eqnarray}
Set $\bm{v} = \left( \begin{array}{c} u +\frac{r}{2} k_{11} k_{12} \\ v + \frac{r}{2}  k_{21} k_{22} \end{array}\right)$. By (\ref{eqn:SPlus_Bihol_2}) and (\ref{eqn:GPlus_Bihol_GammaR_diffeo}),  there exist an isomorphism $\rho:G^{+}_{N',p',q',r'} \to {G^{+}_{N,p,q,r}}$ defined in Lemma \ref{lem:Gamma_R}. It suffices to show that $\varphi \circ g'_0 \circ \varphi^{-1} = g_0 g_1^{l_1} g_2^{l_2}$ and $\varphi \circ g' \circ \varphi^{-1} = \rho(g')$ for any $g' \in \Gamma
_{G^{+}_{N',p',q',r'}}$. Note that
\begin{eqnarray*}
g_0 g_1^{l_1} g_2^{l_2}(w,z) = g_0 (\mu(g_1 g_2) (w,z)).
\end{eqnarray*}
Hence we obtain $\varphi \circ g'_0 \circ \varphi^{-1} = g_0 g_1^{l_1} g_2^{l_2}$ by (\ref{eqn:SPlus_Bihol_6}) and (\ref{eq:Splus_biHolE}). \\
\indent Now, let $\Gamma'_r$ be a copy of $\Gamma_r$. We define the homomorphism $\mu': \Gamma_{G^{+}_{N',p',q',r'}} \to \Gamma'_r$ and the action of $\Gamma'_r$ on $\HC \times \C$ similarly to (\ref{eqn:Gammar_HCAction}) and (\ref{eqn:GammaR_Action}) respectively. 
Then, $\varphi \circ g' \circ \varphi^{-1} = \rho(g')$ is equivalent to $\varphi (\mu' (g') (w,z)) = \mu(\rho(g')) \varphi (w,z)$.  Thus, it is sufficient to prove that $\varphi \circ (\zeta , y) = (\zeta K, \zeta \bm{v} + \det K y) \circ \varphi$ for any $(\zeta, y) \in \Gamma'_r$. This is equivalent to (\ref{eqn:SPlus_Bihol_4}) and
\begin{eqnarray*}
\frac{c}{\alpha -1} (\eta \bm{b}) \bm{a'} + f \left( \zeta \bm{c'}-  \frac{\theta'}{r} y + \frac{1}{2} (\zeta \bm{a'} ) (\zeta \bm{b'} ) \right) & = &  - \frac{\alpha (\eta \bm{a} ) }{\alpha -1} (\zeta K \bm{b}) + \zeta K \bm{c} - \frac{\theta}{r} ((\det K ) y + \bm{v} ) \nonumber  \\ &+& \frac{1}{2} ( \zeta K \bm{a}) (\zeta K \bm{b}) \label{SPlus:mouiikagenishite},
\end{eqnarray*}
for any $(\zeta, y)$.
By (\ref{eqn:SPlus_Bihol_4}), we obtain $\displaystyle \frac{f}{2} (\zeta \bm{a'} ) (\zeta \bm{b'}) = \displaystyle \frac{1}{2} (\zeta K \bm{a} ) (\zeta K \bm{b})$ and $\displaystyle \frac{f \theta'}{r} =  \displaystyle \frac{\theta}{r} (\det K)$. It remains to show that
\begin{eqnarray}
\label{eqn:GammaR_Bihol_Saigo_1}
K \bm{c} - f \bm{c'} = \frac{c  (\eta \bm{b})}{\alpha -1} \bm{a'}  + \frac{\alpha (\eta \bm{a})}{\alpha -1} (K \bm{b}) + \frac{\theta}{r} \bm{v}.
\end{eqnarray}
By (\ref{SPlus_SPEQ}) and (\ref{eqn:GPlus_Bihol_GammaR_diffeo}), (\ref{eqn:GammaR_Bihol_Saigo_1}) is equivalent to
\begin{eqnarray*}
\theta K \left( \begin{array}{c} l_2 \\ - l_1 \end{array} \right)  =  c (\eta \bm{b}) \bm{a'} - \frac{f}{c} (\eta \bm{a}) \bm{b'}.
\end{eqnarray*}
Hence it suffices to show that $Z (l_2 , -l_1)^{T} = 0$ for
\begin{eqnarray*}
Z = \theta K - c \bm{a'} (b_2,-b_1) + \frac{f}{c} \bm{b'} (a_2,-a_1).
\end{eqnarray*}
This follows from $Z \bm{a} = Z \bm{b} = 0$. \qed \\

%
%

%
%
%
%
%
%
%
%
%
%
%
%
%
%
%
%
%

%
%
%
%
%
%
%
%
%
%
%
\begin{Lemma}
\label{Lem:SplusBiholCol}
Assume that 
\begin{eqnarray*}
r' & = & r,  \\
KNK^{-1} & = &  N',\\
\det K \left( \begin{array}{c} p' \\ q' \end{array} \right)   - K \left(  \begin{array}{c} p \\ q \end{array} \right)  & \in & r \Z^2 + (N' - I) \Z^2.   \\
\end{eqnarray*}
for some $K \in GL(2,\Z)$. Fix any $t' \in \C$ and eigenvectors $\bm{a'},\bm{b'}$ of $N'$ corresponding to $\alpha,1/\alpha$ respectively. Then, $S^{+}_{N,p,q,r,t,\bm{a},\bm{b}}$ is deformation equivalent to $S^{+}_{N',p',q',r',t',\bm{a'},\bm{b'}}$ or $S^{+}_{N',p',q',r',t',-\bm{a'},\bm{b'}}$.
\end{Lemma}
\paragraph{Proof}
We use the isomorphism $\varphi:G^{+}_{N',p',q',r'} \to {G^{+}_{N,p,q,r}}$ given in the proof of Lemma \ref{thm:SPlus_diffeo}. By Lemma \ref{lem:Gamma_R}, the lift of $\varphi|\Gamma_{G^{+}_{N',p',q',r'}}$ satisfies (\ref{eqn:GPlus_Bihol_GammaR_diffeo}). Choose $\hat{t} \in \C$ and eigenvectors $\hat{\bm{a}},\hat{\bm{b}}$ which satisfy (\ref{eqn:SPlus_Bihol_4}) and (\ref{eqn:SPlus_Bihol_6}). Then $S^{+}_{N',p',q',r',\hat{t},\hat{\bm{a}},\hat{\bm{b}}}$ is biholomorphic to $S^{+}_{N,p,q,r,t,\bm{a},\bm{b}}$. Note that the choice of $t$ does not affect the deformation equivalence class. Furthermore, if $c>0$ and $f \in \R \setminus \{ 0 \}$, then $S^{+}_{N,p,q,r,0,\bm{a},\bm{b}}$  and $S^{+}_{N,p,q,r,0,c\bm{a},f\bm{b}}$ are biholomorphic by the map
\begin{eqnarray*}
 S^{+}_{N,p,q,r,0,\bm{a},\bm{b}} \to S^{+}_{N,p,q,r,0,c\bm{a},f\bm{b}}: (w,z) \mapsto (cw,fz).
\end{eqnarray*}
This proves the lemma. \qed \\

Let $S^{+}_{N_i,p_i,q_i,r_i,t_i,\bm{a_i},\bm{b_i}}$ be Inoue surfaces homotopy equivalent to $S^{+}_{N,p,q,r,t,\bm{a},\bm{b}}$, where $i=1,2$. For each $i$, there exists $K_i \in GL(2, \Z)$ and $\delta_i,\epsilon_i \in \{1,-1 \}$ which satisfies (\ref{eqn:SPlus_diffeo1}) - (\ref{eqn:SPlus_diffeo3}) of Lemma \ref{thm:SPlus_diffeo}. By Lemma \ref{Lem:SplusBiholCol}, we can easily show the following corollary.
\begin{Colr}
\label{SPlus_Colr}
If $\det K_1 = K_2, \delta_1 = \delta_2$ and $\epsilon_1 = \epsilon_2$, then $S^{+}_{N_2,p_2,q_2,r_2,t_2,\bm{a_2},\bm{b_2}}$ is deformation equivalent either to  
$S^{+}_{N_1,p_1,q_1,r_1,t_1,\bm{a_1},\bm{b_1}}$ or to $S^{+}_{N_1,p_1,q_1,r_1,t_1,-\bm{a_1},\bm{b_1}}$.
\end{Colr}
\begin{prop}
The number of deformation equivalence classes of Inoue surfaces of type $S^{+}$ homotopy equivalent to $S^{+}_{N,p,q,r,t,\bm{a},\bm{b}}$ is at most 16.
\end{prop}
\paragraph{Proof} We first construct 16 surfaces which are homotopy equivalent to $S^{+}_{N,p,q,r,t,\bm{a},\bm{b}}$. Let $J = \displaystyle \left( \begin{array}{cc} 0 & 1 \\ 1 & 0 \end{array} \right)$. For every $(d,\delta,\epsilon) \in \{ 1,-1 \}^3$, we choose integers $P(d , \delta, \epsilon ) $ and $Q(d , \delta, \epsilon )$ satisfying
\begin{eqnarray*}
\delta \left( \begin{array}{c} P(d , \delta, \epsilon ) \\ Q(d , \delta, \epsilon ) \end{array} \right)   - \epsilon J^d \left(  \begin{array}{c} p \\ q \end{array} \right)  & \in & r \Z^2 + (J^{\frac{-d+1}{2}} N^{\epsilon} J^{\frac{d-1}{2}}) \Z^2.
\end{eqnarray*}
Fix eigenvectors $\bm{a}(d,\delta,\epsilon)$ and $\bm{b}(d,\delta,\epsilon)$ of the matrix $J^{\frac{-d+1}{2}} N^{\epsilon} J^{\frac{d-1}{2}}$ corresponding to eigenvalues $\alpha$ and $\frac{1}{\alpha}$ respectively. Set
\begin{eqnarray*}
S^{+}_{(d,\delta,\epsilon,+)} & = & S^{+}_{ J^{\frac{-d+1}{2}} N^{\epsilon} J^{\frac{d-1}{2}}, P(d , \delta, \epsilon ) ,  Q(d , \delta, \epsilon ) , d \delta r,  \bm{a}(d , \delta, \epsilon ), \bm{b}(d , \delta, \epsilon )}, \\
S^{+}_{(d,\delta,\epsilon,-)} & = & S^{+}_{ J^{\frac{-d+1}{2}} N^{\epsilon} J^{\frac{d-1}{2}}, P(d , \delta, \epsilon ) ,  Q(d , \delta, \epsilon ) , d\delta r,  -\bm{a}(d , \delta, \epsilon ), \bm{b}(d , \delta, \epsilon )}.
\end{eqnarray*}
\indent Now, let $S'$ be an Inoue surface of type $S^{+}$ homotopy equivalent to $S^{+}_{N,p,q,r,t,\bm{a},\bm{b}}$. Write $S' = S^{+}_{N',p',q',r',t',\bm{a'},\bm{b'}}$. Then, there exists $K \in GL(2,\Z), \epsilon,\delta \in \{ 1,-1 \}$ which satisfies (\ref{eqn:SPlus_diffeo1}) - (\ref{eqn:SPlus_diffeo3}) . By Corollary \ref{SPlus_Colr}, $S'$ is deformation equivalent either to $S^{+}_{(\det K,\delta,\epsilon,+)}$ or to $S^{+}_{(\det K,\delta,\epsilon,-)}$. This proves the proposition. \qed

\section{Inoue surfaces of type $S^{-}$}
\indent \indent We first recall the definition of Inoue surfaces of type $S^{-}$. Let $N = (n_{ij} ) \in GL(2, \Z) $ be a matrix with eigenvalues $\alpha>1, -1/\alpha$. We choose eigenvectors $\bm{a} = {}^{T} (a_1, a_2) \in \R^2, \bm{b} = {}^{T}(b_1,b_2) \in \R^2$  of $N$ corresponding to $\alpha, -1/\alpha$ respectively. Fix integers $p,q$ and $r$, where $r \neq 0$. Let $\theta = \det( \bm{a}, \bm{b} )$. We define $c_1, c_2 \in \R$ to be the solution of
\begin{eqnarray*}
(N + I) \left( \begin{array}{c} c_1 \\ c_2 \end{array} \right) +  \left( \begin{array}{c} e_1 \\ e_2 \end{array} \right) - \frac{\theta}{r}  \left( \begin{array}{c} p \\ q \end{array} \right) = 0,
\end{eqnarray*}
where $e_i = \frac{1}{2}n_{i1} (n_{i1} -1) a_1 b_1 + \frac{1}{2}n_{i2} (n_{i2} -1) a_2 b_2 + n_{i1} n_{i2} b_1 a_2$ for $i = 1,2$. 
Let $\GMinus$ be the subgroup of $\Aut( \HC \times \C)$ generated by 
\begin{eqnarray}
g_0: (w,z) & \mapsto & (\alpha w, z), \label{eqn:SMinus_Gen1}\\
g_i : (w,z) & \mapsto & (z + a_i, w + b_i z + c_i ), \\
g_3: (w,z) & \mapsto & (z, w - \frac{\theta}{r}), \label{eqn:SMinus_Gen2}
\end{eqnarray}
where $i=1,2$. The quotient surface $S^{-}_{N,p,q,r,\bm{a},\bm{b}}= \HC \times \C / \GMinus$ is called an Inoue surface of type $S^{-}$. \\
\indent By (\ref{eqn:SMinus_Gen1}) - (\ref{eqn:SMinus_Gen2}), we have the following relations:
\begin{gather}
g_0 g_1 g_0^{-1} = g_1^{n_{11}} g_2^{n_{12}} g_3^p,
g_0 g_2 g_0^{-1} = g_1^{n_{21}} g_2^{n_{22}} g_3^q, g_0 g_3 g_0^{-1} = g_3^{-1},   \label{eqn:SMinus_Rel1}\\
g_1 g_2 g_1^{-1} g_2^{-1} = g_3^r, g_i g_3 = g_3 g_i,  \label{eqn:SMinus_Rel2}
\end{gather}
where $i=1,2$. As in section 2, we can show that any element $g \in G^{-}_{N,p,q,r}$ can be written uniquely as a product $g = \prod_{i=0}^{3} g_i^{n_i}$, where $n_0,n_1,n_2,n_3$ are integers. Thus, $\GMinus$   has a presentation in terms of the generators $g_0,g_1,g_2,g_3$, with relations (\ref{eqn:SMinus_Rel1}) and (\ref{eqn:SMinus_Rel2}). In this case, the center of $\GMinus$ is trivial. \\
\indent Let $\Gamma_{\GMinus}$ be the subgroup of $\GMinus$ defined in  (\ref{eqn:Gamma_Group}). Since $\Gamma_{\GMinus}$ is isomorphic to $\Gamma_{G^{+}_{N,p,q,r}}$, properties (i), (ii), and (iii) of section 2 holds for $\Gamma_{\GMinus}$ as well. The quotient group $G^{-}_{N,p,q,r}/\Gamma_{G^{-}_{N,p,q,r}}$ is an infinite cyclic group generated by the class of $g_0$.
\begin{Lemma}
\label{thm:SMinus_diffeo}
Inoue surfaces $S^{-}_{N,p,q,r,\bm{a},\bm{b}}$ and $S^{-}_{N',p',q',r',\bm{a'},\bm{b'}}$ constructed above are homotopy equivalent if and only if there exists $K \in GL(2,\Z) , \delta \in \{ 1,-1 \}$ such that
\begin{gather}
r'  =  \delta \det K r, \label{eqn:SMinus_diffeo1}\\
N'  =  KN K^{-1}, \label{eqn:SMinus_diffeo2} \\
\delta \left( \begin{array}{c}
p' \\
q'
\end{array} \right) -  K \left(
\begin{array}{c}
p \\
q
\end{array}
 \right)  \in  r \Z^2 + (N' + I) \Z^2. \label{eqn:SMinus_diffeo3}
\end{gather}
\end{Lemma}
\paragraph{Proof} The proof is similar to Lemma \ref{thm:SPlus_diffeo}. Assume that $S^{-}_{N,p,q,r,\bm{a},\bm{b}}$ and $S^{-}_{N',p',q',r',\bm{a'},\bm{b'}}$ are homotopy equivalent. Then, there exists an isomorphism $\varphi :G^{-}_{N',p',q',r'} \to \GMinus$. The properties above implies that there exists $\epsilon, \delta \in \{-1,1 \}$ and integers $k_{ij} (i=1,2, j= 1,2,3)$ such that
\begin{eqnarray*}
\varphi (g'_{0}) &= & g_0^{\epsilon} g_1^{k_{01}} g_2^{k_{02}} g_3^{k_{03}}, \label{eqn:SMinusDiffeoGroup1} \\
\varphi (g'_{i}) & = &  g_1^{k_{i1}} g_2^{k_{i2}} g_3^{k_{i3}}. \label{eqn:SMinusDiffeoGroup2}
\end{eqnarray*}
Suppose $\epsilon = 1$. We have (\ref{eqn:SMinus_diffeo1}), (\ref{eqn:SMinus_diffeo2}) and (\ref{eqn:SMinus_diffeo3}) by a similar calculation in Lemma \ref{thm:SPlus_diffeo}. If $\epsilon = -1$, it follows similarly that $KN^{-1}K^{-1} = N'$. However this is a contradiction, since $N'$ must have an eigenvalue $0 <\alpha <1$. \\
\indent Next we show the converse. We construct an isomorphism $\varphi :G^{-}_{N',p',q',r'} \to \GMinus$ similarly to Lemma \ref{thm:SPlus_diffeo}. Hence, $S^{-}_{N,p,q,r,\bm{a},\bm{b}}$ and $S^{-}_{N',p',q',r',\bm{a'},\bm{b'}}$ are diffeomorphic to each other since they are both solvmanifolds \cite{Hasegawa_Solvmanifold}. \qed \\

\begin{Lemma} 
\label{lem:SMinus_BiholConstruction}
Assume that $r = r'$ and that there exists $K \in GL(2,\Z),c>0, u,v \in \Z$ such that 
\begin{gather*}
\eta \in \Z^2,KNK^{-1}  =  N', \label{eqn:SminusBihol_Cond2}
c \bm{a'}  =  K \bm{a}, f \bm{b'}  =  cK \bm{b} \nonumber,
\end{gather*}
where $f = \displaystyle \frac{\det K \theta}{\theta'}$ and $\eta = (l_1,l_2) \in \Z^2$ is a solution of
\begin{eqnarray*}
 (N' + I ) \left( \begin{array}{c} u +\frac{r}{2} k_{11} k_{12} \\ v + \frac{r}{2}  k_{21} k_{22} \end{array}\right)   =   K \bm{p} - (\det K) \bm{p}' - Kr \left( \begin{array}{c} -l_2 \\ l_1 \end{array} \right)  \label{eqn:SminusBihol_Cond3}.
 \end{eqnarray*}
Then $S^{-}_{N,p,q,r,\bm{a},\bm{b}}$ and $S^{-}_{N',p',q',r',\bm{a'},\bm{b'}}$ are biholomorphic.
\end{Lemma}
To prove this lemma, we apply $\Gamma_r$ for Inoue surfaces of type $S^{-}$. By replacing (8.3) of \cite{Surjective_Isom} by $\left( \begin{array}{ccc} 1 & 0 & 0 \\ 0 & \alpha & 0 \\ 0 & 0 & -1 \end{array} \right)$,  equation (8.4) of \cite{Surjective_Isom} becomes
\begin{eqnarray}
\label{eqn:GammaR_C_to_P_SMinus}
(N + I) \bm{c} = \frac{\theta}{r} \bm{p},
\end{eqnarray} 
where $\bm{c}$ and $\bm{p}$ are column vectors defined in section 2. We define the homomorphism $\mu: \Gamma_{\GMinus} \to \Gamma_r$ by (\ref{eqn:Gamma_to_Gammar}), and the action $\Gamma_r$ on $\HC \times \C$ by  (\ref{eqn:Gammar_HCAction}). We can show the following lemma similarly to Lemma \ref{lem:Gamma_R}:
\begin{Lemma}
Assume that $r = r'$. Let $\rho:{G^{-}_{N',p',q',r'}} \to {G^{-}_{N,p,q,r}}$ be an isomorphism such that $\rho(g'_0) = g_0 g_1^{l_1} g_2^{l_2}$ for some integers $l_1,l_2$. Denote by $(K , \bm{v})\in GL(2,\Z) \times \Z[r/2]^2$, the lift of $\rho|{\Gamma_{{G^{-}_{N',p',q',r'}}}}$. Then $(K , \bm{v})$ satisfies the following condtions:
\begin{gather}
KNK^{-1} = N', \nbsp (N' + I ) \bm{v}  =  K \bm{p} - (\det K) \bm{p}' - Kr \left( \begin{array}{c} -l_2 \\ l_1 \end{array} \right), \label{SMinus_GammaR_equiv2} \\
v_i - (r/2) k_{i1} k_{i2}  \in \Z,  \nbsp \nbsp \nbsp \nbsp \label{Gamma_r_isom22}
\end{gather} 
where $\bm{v} = {}^T (v_1 , v_2) $. Conversely, assume that $(K , \bm{v}) \in GL(2,\Z) \times \Z[r/2]^2$ satisfies (\ref{SMinus_GammaR_equiv2}) and (\ref{Gamma_r_isom22}) for some integers $l_1$ and $l_2$. Let $\varphi:\Gamma_{G^{-}_{N',p',q',r'}} \to \Gamma_{G^{-}_{N,p,q,r}}$ be an isomorphism whose lift corresponds to $(K , \bm{v})$. Then, the map
\begin{eqnarray*}
 \rho:G^{+}_{N',p',q',r'} \to G^{+}_{N,p,q,r}: {g'_0}^{l_0} \gamma \mapsto (g_0 g_1^{l_1} g_2^{l_2})^{l_0} \varphi (\gamma),
 \end{eqnarray*}
 is an isomorphism, where $\gamma \in \Gamma_{G^{-}_{N',p',q',r'}}$.
\end{Lemma}
\paragraph{Proof of Lemma \ref{lem:SMinus_BiholConstruction}}
Let $\varphi:\HC \times \C \to \HC \times \C: (w,z) \mapsto (cw + d, ew + fz + g)$, where
\begin{gather*}
d = \frac{\alpha}{1 - \alpha} \eta \bm{a}, \nbsp e = c \frac{\eta \bm{b}}{\alpha +1}, \\g = \frac{1}{2} ( -d \eta \bm{b} - \eta \bm{c} + \frac{\theta}{2}  l_1 l_2) + \frac{1}{4} (\eta \bm{a})(\eta \bm{b}).
\end{gather*}
Set $\bm{v} = \left( \begin{array}{c} u +\frac{r}{2} k_{11} k_{12} \\ v + \frac{r}{2}  k_{21} k_{22} \end{array}\right)$. It suffices to show that $\varphi \circ g'_0 \circ \varphi^{-1} = g_0 g_1^{l_1} g_2^{l_2}$ and $\varphi \circ (\zeta , y) = (\zeta K, \zeta \bm{v} + \det K y) \circ \varphi$ for any $(\zeta, y)  \in \Gamma_r$.

This can be proved similarly to Lemma \ref{lem:SPlusBihol} by using (\ref{eqn:GammaR_C_to_P_SMinus}) and (\ref{SMinus_GammaR_equiv2}).

\qed \\

\begin{Colr}
Assume that there exists $K \in GL(2, \Z)$ such that
\begin{gather*}
r'  =  r, KNK^{-1}  =   N',\\
\det K \left( \begin{array}{c} p' \\ q' \end{array} \right)   - K \left(  \begin{array}{c} p \\ q \end{array} \right)   \in  r \Z^2 + (N' +I) \Z^2.
\end{gather*}
Fix any eigenvectors $\bm{a'},\bm{b'}$ of $N'$ corresponding to $\alpha, -1/\alpha$ respectively. Then, $S^{-}_{N,p,q,r,\bm{a},\bm{b}}$ is deformation equivalent either to $S^{-}_{N',p',q',r',\bm{a'},\bm{b'}}$ or $S^{-}_{N',p',q',r',-\bm{a'},\bm{b'}}$.
\end{Colr}
\begin{prop} \label{prop:SMinusDN}The number of deformation equivalence classes of Inoue Surfaces of type $S^{-}$ homotopy equivalent to $S^{-}_{N,p,q,r,\bm{a},\bm{b}}$ is less than or equal to 8.
\end{prop}
The statements above can be shown similarly to section 2.
\section{Deformation equivalence classes of Inoue surfaces of type $S^0$,$S^{+}$ or $S^{-}$}
\indent \indent In this section, we first review the classification of compact surfaces with $b_1=1, b_2 =0$. Then, we show that the number of deformation equivalence classes of surfaces homotopy equivalent to an Inoue surface of type $S^0$,$S^{+}$, or $S^{-}$ is at most 16. 
\begin{prop}
\label{prop:classification}
Let $S$ be a compact surface with $b_1(S) = 1$ and $b_2(S)  = 0$. Then $S$ is either an elliptic surface, a Hopf surface, or an Inoue surface of type $S^0$,$S^{+}$ or $S^{-}$. 
\end{prop}
\paragraph{Proof} Let $\kappa(S)$ be the Kodaira dimension of $S$. Since $b_1(S)$ is odd, $\kappa(S)$ is either $- \infty,0$ or $1$. If $\kappa(S)= 1$, then $S$ is elliptic \cite{VanDeVan}. If $\kappa(S) = 0$, $S$ must be a secondary Kodaira Surface, which are elliptic \cite{VanDeVan}. By Bogomolov \cite{Bogomolov} and Teleman \cite{Teleman}, if $\kappa(S) = -\infty$, $S$ is either a Hopf surface or an Inoue surface of type $S^0$,$S^{+}$, or $S^{-}$. \qed
\\

By \cite{Inoue_Original} and the results of sections 1 and 2, we see that $G_M,G^{+}_{N,p,q,r},G^{-}_{N,p,q,r}$ has the following properties:
\begin{center}
\begin{tabular}{|c|c|c|c|}
\hline
$G_M$  & trivial center  & $\Gamma_{G_M}$& abelian\\
\hline
$G^{+}_{N,p,q,r}$  & infinite-cyclic center & $\Gamma_{G^{+}_{N,p,q,r}}$& non-abelian \\
\hline
$G^{-}_{N,p,q,r}$ & trivial center & $\Gamma_{G^{-}_{N,p,q,r}}$& non-abelian  \\
\hline
\end{tabular}
\\ 
\vspace{0.5cm}
Table 4.2
\end{center}
Let $S$ be an Inoue surface of type $S^0$,$S^{+}$, or $S^{-}$, and $S'$ be an elliptic surface with $b_1(S') = 1, b_2(S') = 0$. By the proof of [3  Chapter II Theorem 7.16], if $\pi_1(S')$ is non-abelian, then the center of $\pi_1(S')$ is isomorphic to either $\Z \oplus \Z / n \Z  (n \geq 1)$ or $\Z^2$. Hence, $\pi_1(S)$ is not isomorphic to $\pi_1(S')$. Recall that a Hopf surface is by definition, a surface with universal cover $\C^2 \setminus \{ 0 \}.$  Thus, $S$ cannot deform to elliptic surfaces or to Hopf surfaces. Furthermore, Proposition \ref{prop:classification} and Table 4.2 implies that any surface deformation equivalent to $S$ must deform through Inoue surfaces of the same type. Therefore, the number of deformation equivalence classes of surfaces homotopy equivalent to an Inoue surface with $b_2 =0$ is at most 16.
\section{Elliptic Surfaces and Hopf Surfaces}
\indent \indent To complete the proof of our main theorem, we must consider the case when $S$ is a Hopf surface or an elliptic surface with $b_1 (S) = 1$ and $b_2 (S)= 0$. According to [3. Chapter I Lemma 7.20], if $S$ is a Hopf surface with a non-abelian fundamental group, then $S$ is elliptic. It is also known that if $S$ is an elliptic surface whose fundamental group is abelian, then $S$ is a Hopf Surface [3. Chapter II Proposition 7.5]. Hence, we may assume that $S$ is either an elliptic surface whose fundamental group is non-abelian, or a Hopf surface with an abelian fundamental group.
\begin{prop}
Let $S$ be an elliptic surface with $b_1(S) = 1, b_2(S) = 0$. Assume that $\pi_1(S)$ is non-abelian. Then the number of deformation equivalence classes of surfaces homotopy equivalent to $S$ is at most two.
\end{prop}
\paragraph{Proof}
According to \cite{Kodaira_Basic2}, $S$ is obtained by performing logarithmic transformation finitely times over $\CP^1 \times E$, where $E$ is a general elliptic curve. Since the base $\CP^1$ is simply connected, $S$ is an elliptic surface with cyclic monodromy. Let $S'$ be a surface homotopy equivalent to $S$. By [3.   Chapter II Corollary 7.17],  $S'$ is deformation equivalent  either to $S$ or to $S^{conj}$ where $S^{conj}$ is the conjugate complex manifold of $S$ (See [3. Chapter II Definition 7.13] for the definition of $S^{conj}$). \qed 
\begin{prop}
Let  $S$ be a Hopf surface whose fundamental group is abelian. Then the number of deformation equivalence classes of surfaces homotopy equivalent to $S$ is finite.
\end{prop}
\paragraph{Proof}
We refer to [3. Chapter I Section 1.7.6] for details. If $\pi_1(S) \cong \Z$, the number of deformation equivalence classes of surfaces homotopy equivalent to $S$ is one since Hopf surfaces with an infinite-cyclic fundamental group are all deformation equivalent. If $\pi_1(S) \cong \Z \oplus \Z / n \Z (n>1)$, $S$ is diffeomorphic  to $S_1 \times L(n,q)$ for some $q \in (\Z / n \Z)^{*}$, where $L(n,q)$ is a lens space. Since the choice of $q$ is finite, the number of deformation equivalence classes of surfaces homotopy equivalent to $S$ must be finite. 
\qed
\paragraph{Remark} By the discussion of \cite{Friedman}, the number of deformation equivalence classes of surfaces diffeomorphic to $S$ is at most two.\\

This completes the proof of our main theorem.
\providecommand{\bysame}{\leavevmode\hbox to3em{\hrulefill}\thinspace}
\providecommand{\MR}{\relax\ifhmode\unskip\space\fi MR }
\providecommand{\MRhref}[2]{%
  \href{http://www.ams.org/mathscinet-getitem?mr=#1}{#2}
}
\providecommand{\href}[2]{#2}

\noindent Present Address: \\
Department of Mathematics, \\
Keio University,\\
Hiyoshi, Kohoku-Ku, Yokohama, 223-8532 Japan. \\
{\it E-mail address}: ua592825@keio.jp


\begin{thebibliography}{10}

\bibitem{VanDeVan}
W.~Barth, K.~Hulek, C.A.M. Peters, and A.~Van~De Van, \emph{Compact complex
  surfaces}, Springer, 2004.

\bibitem{Bogomolov}
F.~A. Bogomolov, \emph{Classification of surfaces of class ${VII}_0$ with $b_2
  =0$}, Seriya Matematicheskaya 10 (1976).

\bibitem{Friedman}
R.~Friedman and J.W. Morgan, \emph{Smooth 4-manifolds and complex surfaces},
  Springer, 1994.

\bibitem{Surjective_Isom}
Y.~Fujimoto and N.~Nakayama, \emph{Compact complex surfaces admitting
  non-trivial surjective endomorphisms}, Tohoku Mathematical Journal (2005).

\bibitem{Hasegawa_Solvmanifold}
K.~Hasegawa, \emph{Complex and {K}$\ddot{a}$hler structures on compact
  solvmanifolds}, Journal of Symplectic Geometry Volume 3 (2005).

\bibitem{Inoue_Original}
M.~Inoue, \emph{On surfaces of {c}lass ${V}{I}{I}_0$}, Inventiones math (1974).

\bibitem{Inoue_SM}
\bysame, \emph{An example of an analytic surface}, Sugaku 27, (1975), 358--364
  (in Japanese).

\bibitem{Kato}
M.~Kato, \emph{Compact complex manifolds containing "global" spherical shells.
  {I}}, Proceedings of the International Symposium on Algebraic Geometry
  (1978).

\bibitem{Kodaira_Basic2}
K.~Kodaira, \emph{On the structure of compact complex analytic surfaces,
  {II}.}, American Journal of Mathematics 86 (1966).

\bibitem{Teleman}
A.~Teleman, \emph{Projectively flat surfaces and {B}ogomolovfs theorem on
  {C}lass ${VII}_0$}, International Journal of Mathematics (1994).

\bibitem{Teleman_arxiv}
\bysame, \emph{Gauge theoretical methods in the classification of
  non-kaehlerian surfaces}, Arxiv preprint (2008).

\bibitem{Yau}
S.~T. Yau, \emph{Calabifs conjecture and some new results in algebraic
  geometry}, Proceedings of the National Academy of Sciences of the United
  States of America 74 (1977).

\end{thebibliography}
\end{document}